\documentclass[12pt]{amsart}
\usepackage[top=1in, bottom=1in, left=1in, right=1in]{geometry}

\usepackage{amssymb,mathrsfs,amsmath}

\newtheorem{theorem}{Theorem}[section]

\theoremstyle{remark}

\newcommand{\bsp}{\begin{split}}
\newcommand{\esp}{\end{split}}
\newcommand{\be}{\begin{equation}}
\newcommand{\ee}{\end{equation}}
\newcommand{\bes}{\begin{equation*}}
\newcommand{\ees}{\end{equation*}}
\newcommand{\bv}\boldsymbol{}

\numberwithin{equation}{section}

\renewcommand{\pmod}[1]{~({\rm mod}\,#1)}

\begin{document}

\title{Big biases amongst products of two primes}
\author{David Dummit}
\address{Department of Mathematics and Statistics\\
University of Vermont\\
Burlington, VT 05401 \\
USA}
\email{{\tt dummit@math.uvm.edu}}
\author{Andrew Granville}
\address{Centre de recherches math\'ematiques\\
Universit\'e de Montr\'eal\\
CP 6128 succ. Centre-Ville\\
Montr\'eal, QC H3C 3J7\\
Canada}
\email{{\tt andrew@dms.umontreal.ca}}
\author{Hershy Kisilevsky}
\address{Department of Mathematics and Statistics\\
Sir George Williams Campus\\
Concordia University\\
Montreal, QC H3G 1M8 \\
Canada}
\email{{\tt kisilev@mathstat.concordia.ca  }}

\subjclass[2010]{11L20, 11N25, 11N13}
\keywords{Primes in arithmetic progressions, Prime races}

\date{\today}

\begin{abstract} We show that substantially more than a quarter of the odd integers of the form $pq$ up to $x$, with $p,q$ both prime, satisfy $p\equiv q\equiv 3 \pmod 4$.
\end{abstract}

\maketitle


\newcommand{\eS}{\mathcal S}
\newcommand{\eL}{\mathcal L} 
\newcommand{\Lam}{\Lambda}
 
\section{Introduction}  There are roughly equal quantities of odd integers $n$ that are the product of two primes, $p$ and $q$, in the two arithmetic progressions $1\pmod 4$ and $3\pmod 4$. Indeed the counts differ by no more than $x^{1/2+o(1)}$
(assuming  the Riemann Hypothesis for $L(1,(-4/.))$; see \cite{FS} for a detailed analysis).
One might guess that these integers are further evenly  split  amongst those with $(p\pmod 4,q\pmod 4)=(1,1), (1,-1), (-1,1)$ or $(-1,-1)$, but recent calculations reveal a substantial bias towards those $pq\leq x$ with $p\equiv q\equiv 3 \pmod 4$.  Indeed for the ratio
\[
 r(x):=\#\{ pq\leq x: \ p\equiv q\equiv 3 \pmod 4\} \bigg/ \frac 14 \#\{ pq\leq x \} 
\]
we found that 
\[ 
r(1000)\approx 1.347, \ r(10^4)\approx 1.258, \ r(10^5)\approx 1.212, \ r(10^6)\approx 1.183, \ r(10^7)\approx 1.162,
 \]
 a pronounced bias that seems to be converging to $1$ surprisingly slowly. We will show that this is no accident and that there is similarly slow convergence for many such questions:

\begin{theorem} Let $\chi$ be a quadratic character of conductor $d$. For $\eta=-1$ or $1$ we have 
 \[
 \frac{\#\{ pq\leq x: \ \chi(p)=\chi(q)=\eta\}}{\frac 14 \#\{ pq\leq x: \ (pq,d)=1\}} = 
1+\eta \frac{ (\mathcal{L}_\chi +o(1))}{\log\log x}  \ \ \text{where} \ \  \mathcal{L}_\chi := \sum_{ p   }   \frac {\chi(p)}{p }   .
 \]
\end{theorem}

If $\chi=(-4/.)$ then $\mathcal{L}_\chi=-.334\ldots$ so the theorem implies that $r(x)\geq 1+\frac{(1+o(1))}{3(\log\log x-1)}$. If we let $s(x)=1+\frac{1}{3(\log\log x-1)}$ then we have 
\[
s(1000)\approx 1.357, \ s(10^4)\approx 1.273, \ s(10^5)\approx 1.230, \ s(10^6)\approx 1.205, \ s(10^7)\approx 1.187,
\]
 a pretty good fit with the data above. The prime numbers have only been computed up to something like $10^{24}$ so it is barely feasible that one could collect data on this problem up to $10^{50}$ in the foreseeable future. Therefore we would expect this bias to be at least $7\% $ on any data that will be collected this century (as $s(10^{50})\approx 1.07$).

\begin{proof}  For a given quadratic Dirichlet character $\chi$ we will count the number of integers $pq\leq x$ with $\chi(p)=\chi(q)=1$ (and the analogous argument works for $-1$). One can write any such integer $pq\leq x$ with $p\leq q\leq x/p$, so that $p\leq \sqrt{x}$. Hence we wish to determine
\[
 \sum_{\substack{p\leq \sqrt{x} \\ \chi(p)=1}} \sum_{\substack{p\leq q\leq x/p \\ \chi(q)=1}} 1.
\]
The prime number theorem for arithmetic progressions reveals that 
$\sum_{q\leq Q,\ \chi(q)=1} 1 = \frac Q{2\log Q}+O(\frac Q{(\log Q)^2})$, so the above sum equals
\[
  \sum_{ p\leq \sqrt{x}  } \frac{(\chi_0(p)+\chi(p))}2\cdot \frac x{2p\log (x/p)}+O\left(\frac x{p(\log x)^2} + \frac{p}{\log p}\right)
\]
where the implicit constant in the $O(.)$ depends only on the conductor of $\chi$.  This equals
\[
 \frac 14 \sum_{\substack{p\leq \sqrt{x} \\ (p,d)=1}}   \frac x{p\log (x/p)}+ \frac x4  \sum_{ p\leq \sqrt{x}  }   \frac {\chi(p)}{p\log (x/p)}+O\left(\frac x{(\log x)^2} \log\log x\right) .
\]
The difference between the second sum, and the same sum with $\log(x/p)$ replaced by $\log x$, is
\[
 \frac x{4\log x}  \sum_{ p\leq \sqrt{x}  }   \frac {\chi(p)}{p\log (x/p) } \ll \frac x{(\log x)^2}.
\]
using the prime number theorem for arithmetic progression and partial summation. Moreover
\[
  \frac x{4\log x}  \sum_{ p> \sqrt{x}  }   \frac {\chi(p)}{p }\ll \frac x{(\log x)^2}.
\]
Collecting together what we have proved so far yields that $\#\{ pq\leq x: \ \chi(p)=\chi(q)=1\}$
\[
  = \frac 14 \left\{ \#\{ pq\leq x:\ (p,d)=1 \} +
 \frac x{\log x}  \sum_{ p   }   \frac {\chi(p)}{p } +O\left(\frac x{(\log x)^2} \log\log x\right)  \right\}
\]
The first term is well-known to equal $\frac x{\log x}( \log\log x+O(1))$, and so we deduce that 
\[
 \frac{4\#\{ pq\leq x: \ \chi(p)=\chi(q)=1\}}{\#\{ pq\leq x: \ (pq,d)=1\}} = 
1+ \frac{1}{\log\log x} \left(  \sum_{ p   }   \frac {\chi(p)}{p }  +o(1) \right) .
\]
as claimed. 
\end{proof}

We note that 
\[
 \sum_{ p   }   \frac {\chi(p)}{p } =  \sum_{m\geq 1} \frac{\mu(m)}m \log L(m,\chi^m) = 
\log L(1,\chi)+E(\chi),
\]
where  
\[
\sum_p \left( \log \left( 1 -\frac 1p\right) +\frac 1p\right) 
 = -0.315718\ldots \leq E(\chi) \leq \sum_p \left( \log \left( 1 +\frac 1p\right) -\frac 1p\right) 
= -0.18198\ldots
\]
  
\section{Further remarks}

\medskip
$\bullet$\ One deduces from our theorem that $r(x)> 1$ for all $x$ sufficiently large and we conjecture that this is true for all $x\geq 9$.

\medskip
$\bullet$\ We also conjecture that $\mathcal{L}_\chi$ is always non-zero so that there is always such a bias.

\medskip
$\bullet$\ One can calculate the bias in other such questions. For example, we get roughly triple the bias for the proportion of $pq\leq x$ for which
$\left( \frac p5\right)=\left( \frac q5\right)=-1$ out of all   $pq\leq x$ with $p,q\ne 5$ (since $\mathcal{L}_{(./5)}\approx -1.008$). The data
\[ 
r_5(1000)\approx 1.881, \ r_5(10^4)\approx 1.626, \ r_5(10^5)\approx 1.523, \ r_5(10^6)\approx 1.457, \ r_5(10^7)\approx 1.416,
 \]
 confirms this very substantial bias. It would be interesting to find  more extreme examples.

\medskip
$\bullet$\ How large can the bias get if $d\leq x$? It is known \cite{GS} that $L(1,\chi)$ can be as large as $c \log\log d$, and so 
$\mathcal{L}_\chi$ can be as large as $\log\log\log d+O(1)$.  We conjecture that there exists $d\leq x$ for which the bias in our Theorem is as large as 
\[
 1+  \frac{\log\log\log x +O(1)}{\log\log x}  .
\]
Note that this requires proving a uniform version of the Theorem. Our proof assumes that $x$ is allowed to be very large compared to $d$, so does not immediately apply to the problem that we have just stated.

\medskip
$\bullet$\ The same bias can be seen (for much the same reason) in looking at
\[
\sum_{\substack{p\leq x \\ p\equiv 3 \pmod 4}} \frac 1p \ \bigg/  \sum_{\substack{p\leq x \\ p\equiv 1 \pmod 4}} \frac 1p \approx  1+ \frac{2}{3\log\log x}  .
\]
Indeed, by the analogous proof, we have in general
\[
\sum_{\substack{p\leq x \\ \chi(p)=1}} \frac 1p \ \bigg/  \sum_{\substack{p\leq x \\ \chi(p)=-1}} \frac 1p = 1+2 \frac{ (\mathcal{L}_\chi +o(1))}{\log\log x}  .
\]
We therefore see a bias in the distribution of primes in arithmetic progressions, where each prime $p$ is weighted by $1/p$, as a consequence of the sign of $\mathcal{L}_\chi$. This effect is much more pronounced than in the traditional prime race problem where the same comparison is made, though with each prime weighted by $1$.  The bias here is determined by the distribution of values of $\chi(p)$, whereas the prime race bias is determined by the values of $\chi(p^2)=1$, so they appear to be independent phenomena. However one might guess that both biases are sensitive to low lying zeros of $L(s,\chi)$. This probably deserves further investigation, to determine whether there are any correlations between the two biases.

\medskip
$\bullet$\ 
One can show the following for $k$ prime factors, by similar methods:
\[
 \frac{\#\{ p_1\ldots p_k\leq x: \ \text{each} \ \chi(p_j)=\eta\}}{\frac 1{2^k} \#\{  p_1\ldots p_k\leq x: \ \text{each} \ (p_j,d)=1\}} = 
1+\eta \frac{ ((k-1)\mathcal{L}_\chi +o(1))}{\log\log x}    .
 \]
It would be interesting to understand this when $k$ gets large, particularly when $k\sim \log\log x$, the typical number of prime factors of an integer $\leq x$.  It seems likely that the factor on the right-side should grow like
\[
 c_{\chi,k} \left( 1+\eta \frac{ ( \mathcal{L}_\chi +o(1))}{\log\log x}  \right)^{k-1} ,
\]
but we do not know what $c_{\chi,k}$ would look like.

\medskip
$\bullet$\ More generally, if $\chi_1,\ldots ,\chi_k$ are quadratic characters (with $\chi_j$ of conductor $d_j$), and each $\eta_j=-1$ or $1$ then 
\[
 \frac{\#\{ p_1\ldots p_k\leq x: \ \chi_j(p_j)=\eta_j \ \text{for each} \ j\}}{\frac 1{2^k} \#\{  p_1\ldots p_k\leq x: \ \text{each} \ (p_j,d_j)=1\}} = 
1+  \frac{ ((k-1)c(\vec{\chi},\vec{\eta}) +o(1))}{\log\log x}  .
 \]
where 
\[ c(\vec{\chi},\vec{\eta}):= \frac 1k \sum_{j=1}^k \eta_j \mathcal{L}_{\chi_j}.\]
In particular this type of bias does not appear when  $k=2,\ \chi_1=\chi_2$ and $\eta_1+\eta_2=0$.  
Can one prove that $c(\vec{\chi},\vec{\eta})$ can only be $0$ for such trivial reasons? That is, is $c(\vec{\chi},\vec{\eta})=0$ if and only if
$\sum_{j:\ \chi_j=\chi} \eta_j=0$  for every character $\chi\in \vec{\chi}$?

\medskip
$\bullet$\ Given arithmetic progressions $a \pmod m$ and $b\pmod n$, one can surely prove that there exists $\beta=\beta( a \pmod m, b\pmod n)$ such that 
\[
 \frac{\#\{ pq\leq x: \ p\equiv a \pmod m,\  q\equiv b\pmod n\}}{\frac 1{\phi(m)\phi(n)} \#\{ pq\leq x: \ (p,m)=(q,n)=1\}} = 
1+  \frac{ \beta+o(1)}{\log\log x}  .
 \]
It would be interesting to classify when $\beta( a \pmod m, b\pmod n)$ is non-zero, and to determine situations in which it is large. Or more generally for what subsets $A\subseteq (\mathbb Z/m\mathbb Z)^*$ and 
  $B\subseteq (\mathbb Z/n\mathbb Z)^*$ is there no such bias? We would guess that this would only be the case   if either

(i)\ $A$ and $B$ both contain all congruence classes (that is, every prime not dividing $mn$ can be   represented by both $A$ and $B$); or

 (ii)\  $A\cup B$ is a partition of the integers coprime to $mn$ (that is, every prime not dividing $mn$ is   represented by $A$, or represented by $B$, but not both).

\end{document}